\documentstyle[12pt]{article}

\newtheorem{prop}{Proposition}

\newtheorem{lem}{Lemma}
\newtheorem{thm}{Theorem}

\newtheorem{cor}{Corollary}

\newcommand{\be}{\begin{equation}}
\newcommand{\ee}{\end{equation}}

\def\real{\hbox{\rm\setbox1=\hbox{I}\copy1\kern-.45\wd1 R}}

\def\R{\hbox{\rm\setbox1=\hbox{I}\copy1\kern-.45\wd1 R}}

\def\prob{\hbox{\rm\setbox1=\hbox{I}\copy1\kern-.45\wd1 P}}

\def\Zint{{\mathchoice{\setbox1=\hbox{\sf Z}\copy1\kern-.75\wd1\box1}
         {\setbox1=\hbox{\sf Z}\copy1\kern-.75\wd1\box1}
         {\setbox1=\hbox{\footnotesize\sf Z}\copy1\kern-.75\wd1\box1}
         {\setbox1=\hbox{\footnotesize\sf Z}\copy1\kern-.75\wd1\box1}}}

\def\bx{{\bf x}}

\def\by{{\bf y}}

\newcommand{\deq}{\stackrel{\scriptscriptstyle\triangle}{=}}

\newcommand{\F}{\mbox{$\cal F$}}

\newcommand{\HC}{\mbox{$\cal H$}}

\newcommand{\hh}{\mbox{$\hat{h}$}}

\begin{document}
\bibliographystyle{plain}

\thispagestyle{empty}
\setcounter{page}{0}

\vspace {2cm}

{\Large  A. B. Nobel, G. Morvai, S.R. Kulkarni: }

\vspace {1cm}
{ \Huge 
Density estimation from an individual numerical sequence. }

\vspace {1cm}

{\Large IEEE Trans. Inform. Theory  44  (1998),  no. 2, 537--541.
}

\vspace {2cm}

\begin{abstract}
This paper considers estimation of a univariate density 
from an individual numerical sequence.  It is assumed that 
(i) the limiting relative frequencies of the numerical sequence
are governed by an unknown density, and (ii) 
there is a known upper bound for the variation of the density on
an increasing sequence of intervals.
A simple estimation scheme is proposed, and is shown
to be $L_1$ consistent when (i) and (ii) apply.  In addition it is
shown that there is no consistent estimation scheme for the set of 
individual sequences satisfying only condition (i).
\end{abstract}


\pagebreak

\newpage

\section{Introduction}
\label{INT}

Estimation of a
univariate density from a finite data set
is an important problem in theoretical and applied statistics.  
In the most common setting, it is assumed that data are obtained from
a stationary process $X_1, X_2, \ldots$ such that 
\[
\prob\{ X_i \in A \} = \int_A f dx
\ \ \mbox{ for every Borel set } \ \ A \subseteq \real
\]
i.e.\ the common distribution 
of the $X_i$ has density $f$, written $X_i \sim f$.
For each $n \geq 1$ an estimate
$\hat{f}_n$ of $f(\cdot)$ is produced from $X_1, \ldots, X_n$.
The estimates $\{ \hat{f}_n \}$ are said to be 
strongly $L_1$ consistent if $\int |\hat{f}_n - f| dx \to 0$ as 
$n \to \infty$ with probability one.

Common density estimation methods include histogram, 
kernel, nearest neighbor, orthogonal series, wavelet, spline,
and likelihood based procedures.  
For an account of these methods, we refer the interested reader to 
the texts of Devroye and Gy\"orfi \cite{DevGyo85}, 
Silverman \cite{Sil86}, Scott \cite{Scott92}, and
Wand and Jones \cite{WanJon95}.
In establishing consistency and rates of convergence
for estimation procedures like those above, many analyses assume that 
$X_1, X_2, \ldots$ are independent
and identically distributed (i.i.d.), in which
case the distribution of the
process $\{ X_i \}$ is completely specified by the marginal density
$f$ of $X_1$.

Complementing work for independent random variables,
numerous results have also been obtained for 
stationary sequences exhibiting both short and long range 
dependence.
Roussas \cite{Rou67} and Rosenblatt \cite{Ros70} 
studied the consistency and
asymptotic normality of kernel density estimates from Markov
processes.  Similar results, under 
weaker conditions, were obtained by Yakowitz \cite{Yak89}.
Gy\"orfi \cite{Gyo81} showed that there is a simple kernel-based 
procedure $\Phi$ that is strongly $L_2$-consistent for 
every stationary ergodic process $\{X_i\}_{i=-\infty}^{\infty}$
such that (i) the conditional distribution of $X_1$ given 
$\{X_i : i \leq 0\}$ is absolutely continuous with probability one, 
and (ii) the corresponding conditional density $h$ satisfies 
$E \int |h(u)|^2 du < \infty$.
For additional work in this area, see also 
Ahmad \cite{Ahm79}, Castellana and Leadbetter \cite{CasLea86},
Gy\"orfi and Masry \cite{GyoMas90},
Hall and Hart \cite{HalHar90}, and the references contained therein.

With these positive results have come examples showing that
density estimation from strongly dependent processes can be
problematic.  
In a result attributed to Shields, it was shown by 
Gy\"orfi, H\"ardle, Sarda and Vieu \cite{GHSV89} 
that there are histogram density estimates, consistent for every
i.i.d.\ process, that fail for some stationary ergodic process.
Gy\"{o}rfi and Lugosi \cite{GyoLug92} established a similar result
for ordinary kernel estimates.
Extending these results,
Adams and Nobel \cite{AdaNob97} have recently 
shown that there is no density
estimation procedure that 
is consistent for every stationary ergodic process. 

With a view to considering density estimation in a more general
setting, one may eliminate stochastic assumptions.
Here we consider the estimation of an unknown density from
an individual numerical sequence, which need not be the trajectory of
a stationary stochastic process.  
We propose a simple
estimation procedure that is applicable in a purely deterministic
setting.  This deterministic point of view is in line with
recent work on individual sequences in information theory, statistics,
and learning theory 
(cf.\ \cite{Ziv78,MerFedGut92,KulPos95,HauKivWar94}).  
Extending the techniques developed in this paper,
Morvai, Kulkarni, and Nobel \cite{MorKulNob97} consider the problem of
regression estimation from individual sequences.  

In many cases, results based
on deterministic analyses can be applied to individual sample paths in
a stochastic setting.  Theorem \ref{thm1} of this paper
yields a positive result concerning density estimation from ergodic
processes (see Corollary \ref{cor1} below).

\section{The Deterministic Setting}

Let $f: \real \to \real$ be a univariate density function
with associated probability measure $\mu_f(A) = \int_A f(x)dx$. 
An infinite sequence $\bx = (x_1, x_2, \ldots)$ 
of numbers $x_i \in \real$ has {\em limiting density} $f$ if
\be
\label{lrf}
\hat{\mu}_n(A) \, = \, \frac{1}{n} \sum_{i=1}^n I\{x_i \in A\} 
\, \to \, \mu_f(A) 
\ee
for every interval $A \subseteq \real$.
A sequence $\bx$ 
having a limiting density will be called {\em stationary}. 
Let $\Omega(f)$ be the set of
stationary sequences with limiting density $f$.

Note that stationarity concerns the
limiting behavior of relative frequencies, which need not converge
to their corresponding probabilities at any particular rate.
Stationarity says nothing about the mechanism by which 
the individual sequence $\bx$ is produced.
In particular, the limiting relative frequencies 
of a stationary sequence $\bx$ 
are unchanged if one appends to $\bx$ a prefix of any finite length.

The sample paths of ergodic processes
provide one source of stationary sequences.
The next proposition follows easily from 
Birkhoff's ergodic theorem.

\begin{prop}
\label{birk}
If $X_1, X_2, \ldots$ 
are stationary and ergodic with $X_i \sim f$, then
${\bf X} = (X_1, X_2, \ldots) \in \Omega(f)$
with probability one.
\end{prop}

A univariate density estimation scheme is a countable collection $\Phi$ 
of Borel-measurable mappings 
$\phi_n: \real \times \real^n \to \real$, $n \geq 1$.
Thus $\phi_n$ associates every vector 
$(x_1, \ldots, x_n) \in \real^n$ with a 
function $\phi_n(\cdot : x_1, \ldots, x_n)$, 
which is viewed as the estimate of an unknown density 
associated with the sequence $x_1, \ldots, x_n$.
These estimates may take negative values, and they
need not integrate to one.  In particular, no regularity conditions are
imposed on the behavior of $\phi_n$ as a function of its inputs.

A scheme $\Phi$ is $L_1$ consistent for a 
a collection $\Omega$ of stationary sequences
if for each $\bx \in \Omega$,
\[
\int | \phi_n(x : x_1,\ldots,x_n) - f(x) | dx \to 0 \, ,
\]
as $n \to \infty$, where $f$ is the limiting density of $\bx$.
A scheme $\Phi$ is universal if it is $L_1$ consistent for the
set $\Omega^*$ of all stationary sequences.
Note that, for i.i.d.\ data,
a density estimation scheme is called universal if it 
is consistent for every marginal density $f$.
The notion of universality defined above is considerable stronger, as
there are no constraints apart from stationarity
placed on the structure of the individual sequences.
In what follows, when $\bx = x_1, x_2, \ldots$ is fixed, 
$\phi(x : x_1,\ldots,x_n)$ will be denoted by $\phi_n(x)$.

Recall that the total variation of a real-valued
function $h$ defined on an interval $[a,b) \subseteq \real$
is given by
\[
V(h:a,b) = \sup \sum_{i=1}^n |h(t_i) - h(t_{i-1})| \, ,
\]
where the supremum is taken over all finite ordered sequences
$a \leq t_0 < \cdots < t_n < b$.
For each nondecreasing function 
$\alpha: \Zint^+ \rightarrow (0, \infty)$ 
let $\F(\alpha)$ be the set of all densities $f$ on $\real$ such that
$V(f:-i,i) < \alpha(i)$ for $i \geq 1$, and let
\[
\Omega(\alpha) = \bigcup_{f \in {\cal F}(\alpha)} \Omega(f)
\]
be the collection of all those
stationary sequences having limiting densities in $\F(\alpha)$.

\vskip.1in

Given a function $\alpha(\cdot)$ as above,
we propose a simple histogram based procedure that is consistent 
for $\Omega(\alpha)$. 
For each $k \geq 1$ let $\pi_k$ be the partition of $\real$
into dyadic intervals of the form
\[
A_{k,j} = \left[ \frac{j}{2^k} , \frac{j+1}{2^k} \right) 
\ \ \mbox{ with } \ \ j \in \Zint \, ,
\]
and let $\pi_k[x]$ be the unique cell of $\pi_k$ containing 
$x$.  Let $\{b_n\}$ be any sequence of positive integers tending to infinity.
For each sequence of numbers $x_1, \ldots, x_n$ 
and each $k \geq 1$ define histogram density estimates
\be
\label{fkn}
\hh_{n,k}(x) = \frac{2^k}{n} \sum_{i=1}^n I\{ x_i \in \pi_k[x] \} \, .
\ee
Our estimate is selected from among the histograms
$\hat{h}_{n,k}$ by selecting a suitable
value of $k$.
Find the partition index
\be
\label{def2}
k_n = \max\left\{ 1\leq k \leq b_n  : V(\hh_{n,k} : -i,i) < 4 \alpha(i) 
            \ \ \mbox{ for } \ \ 1 \leq i \leq k \right\} \, 
\ee
and define
\be
\label{def3}
\phi_n^*(x : x_1,\ldots,x_n) = \hh_{n,k_n} (x) \, .
\ee
If the conditions defining $k_n$ are not satisfied for any 
$1\leq k \leq b_n$, then set $\phi_n^* \equiv 0$.

\begin{thm}
\label{thm1} 
Let $\alpha: \Zint^+ \rightarrow (0, \infty)$ be a fixed, non-decreasing
function.  The estimation scheme $\Phi^* = \{ \phi_n^* \}$ 
defined by (\ref{fkn})-(\ref{def3}) is $L_1$-consistent
for $\Omega(\alpha)$.
Thus for every stationary sequence $\bx$ 
with limiting density $f \in \F(\alpha)$,
$\int |\phi_n^*(x) - f(x)| dx \to 0$.
\end{thm}

\begin{cor}
\label{cor1}
Let $\alpha(\cdot)$ be fixed and let $\phi_n^*$ be defined by
(\ref{fkn})-(\ref{def3}).
For every stationary ergodic process $\{X_i\}$ 
such that $X_i \sim f$ with $f \in \F(\alpha)$,
\[
\int |\phi_n^*(x : X_1, \ldots, X_n) - f(x)| dx \to 0 
\]
as $n \to \infty$ with probability one.
\end{cor}

\vskip.1in

\noindent
{\bf Example:} Fix $\gamma > 0$, and consider the class of stationary
ergodic processes $\{X_i\}$ such that 
$X_i \sim f$ with $V(f : -\infty,\infty) < 2 \gamma$.  
This class includes, but is not limited to,
processes having uniform, exponential, and normal marginal densities 
with arbitrary means, under the restriction that $Var(X_i)$ 
is greater than $(12 \gamma^2)^{-1}$,  $\gamma^{-2}$, 
and  $(2\pi \gamma^2)^{-1}$, respectively. 
By Corollary \ref{cor1} there is a strongly consistent density estimation
procedure $\Phi^*$ for this class of processes.

\vskip.1in
 
\noindent
{\bf Remark:} The variations used to define $\phi_n^*$ depend on the
cumulative difference between the relative frequencies of adjacent
cells: 
\be
\label{varsum}
V(\hh_{n,k} : -i,i) = 
2^{-k} \sum_{j = -i 2^k}^{i 2^k-2} 
| \hat{\mu}_n(A_{k,j}) - \hat{\mu}_n(A_{k,j+1})| \, .
\ee
To find $\phi_n^*$, put $x_1, \ldots, x_n$ in increasing order,
and then calculate $V(\hat{h}_{n,k}: -i,i)$ 
for each $k = 1, \ldots, b_n$
and each $i = 1, \ldots, k$ by scanning the ordered $x_i$ from left 
to right.  This will require at most $O( n \log n + n b_n)$ operations.

\vskip.1in

In order to apply the procedure $\Phi^*$ described in 
(\ref{fkn})-(\ref{def3}), one must know before seeing $\bx$ that
the variation of its limiting density is
less than a known constant on every interval of the form $[-i,i)$.
The following result shows that this requirement
cannot be materially weakened.

\begin{thm}
\label{thm2}
Let $\F$ be the collection of densities $f$ supported on $[0,1]$
for which $V(f:0,1)$ is finite.
There is no $L_1$ consistent density estimation scheme for 
\[
\Omega = \bigcup_{f \in {\cal F}} \Omega(f) .
\]
In particular, there is no universal density estimation 
scheme for individual sequences.
\end{thm}

If an upper bound on the variance of the unknown density
$f$ were known, the scheme of Theorem \ref{thm1} would 
provide consistent estimates of $f$.

Given any density estimation scheme $\Phi = \{ \phi_n \}$, 
the proof of Theorem 
\ref{thm2} shows how one may construct a stationary 
sequence $\bx$, depending on $\Phi$, 
for which $\phi_n(\cdot)$ fails to converge.
A related argument is used by Adams and Nobel 
\cite{AdaNob97} to show 
that there is no universal density estimation
scheme for stationary ergodic processes.  
As a universal density estimation scheme for individual sequences
would, by virtue of Proposition \ref{birk}, yield a universal scheme
for ergodic processes, their result also implies Theorem \ref{thm2}.

The proof of Theorem \ref{thm1} is given in the next section
after several preliminary results.  The proof of 
Theorem \ref{thm2} is given in Section \ref{PT2}.

\section{Proof of Theorem \protect{\ref{thm1}}}
\label{PT1}

\noindent
{\bf Definition:}
For each partition $\pi$ of $\real$ 
into finite intervals and each $f \in L_1$ define
\[
(f \circ \pi)(x) = \frac{1}{l(\pi[x])} \int_{\pi[x]} f(u) du \, ,
\]                   
where $l(A)$ denotes the length of an interval $A$.
Note that $f \circ \pi$ is piecewise constant on the cells of $\pi$.

\vskip.1in

\begin{lem}
\label{lem2}
Let $\pi_1, \pi_2, \ldots$ be the partitions used to define 
the estimates $\phi^*_n$.
For each pair of integers $k, i \geq 1$,
\[
V(f \circ \pi_k : -i, i) \, \leq \, 3 V(f : -i, i) .
\]
Moreover, if $\bx \in \Omega(f)$ then
\[
\lim_{n \to \infty} V(\hh_{n,k} : -i, i) \, = \,  
V(f \circ \pi_k : -i, i) .
\]
\end{lem}

\noindent
{\bf Proof:} For $f$ non-decreasing it is immediate that 
$V(f \circ \pi_k : -i,i)\le V(f : -i,i)$.  
If $V(f: -i,i) = C < \infty$ then 
$f(x) = u(x) - v(x)$ where $u(\cdot)$ and $v(\cdot)$ are
non-decreasing, $V(u: -i,i) \leq C$ and $V(v: -i,i) \leq 2 C$ 
(cf. Kolmogorov and Fomin \cite{KolFom70}).  
It follows from the definition that 
$f \circ \pi_k = u \circ \pi_k - v \circ \pi_k$, and
since $u$ and $v$ are non-decreasing, so are 
$u \circ \pi_k$ and $v \circ \pi_k$. Therefore
\begin{eqnarray*}
V(f \circ \pi_k : -i, i)&=& V(u \circ \pi_k -v \circ \pi_k: -i, i)\\
&\le&V(u \circ \pi_k : -i, i)+V(v \circ \pi_k : -i, i)\\
&\le& V(u  : -i, i)+V(v  : -i, i)\\
 &\le& 3C
\end{eqnarray*}
as the variation of the sum is less
than the sum of the variations.  
To establish the second claim, note that as $n \to \infty$
\begin{eqnarray*}
V(\hh_{n,k} : -i,i) & = &
2^{-k} \sum_{j = -i 2^k}^{i 2^k - 2} 
| \hat{\mu}_n(A_{k,j}) - \hat{\mu}_n(A_{k,j+1})| \\
& \to &
2^{-k} \sum_{j = -i 2^k}^{i 2^k - 2} 
|\mu_f(A_{k,j}) - \mu_f(A_{k,j+1})| \\
& = & V(f \circ \pi_k: -i,i) \, .
\end{eqnarray*}
$\Box$

\begin{lem}
\label{ktoinf}
Let $\bx \in \Omega(\alpha)$ with limiting density $f \in \F(\alpha)$.
Then the partition index $k_n$ of the density estimate $\phi_n^*$ 
tends to infinity with $n$.
\end{lem}

\noindent
{\bf Proof:}
By Lemma~\ref{lem2}, for arbitrary $K \geq 1$ and 
for all $i = 1,\dots,K$,  
\[\lim_{n \to \infty} V(\hh_{n,K} : -i, i) \, = \,  
V(f \circ \pi_K : -i, i) \le  3 V(f : -i, i)<3\alpha(i) .
\]
Thus by definition of $k_n$, $\liminf_{n\to\infty} k_n \geq K$. $\Box$

\vskip.2in

\noindent
{\bf Proof of Theorem \ref{thm1}:} 
Let $\bx \in \Omega(\alpha)$ be a fixed stationary sequence 
with limiting density $f \in \F(\alpha)$.
For each $n \geq 1$ such that  $k_n \geq 1$ define the error function 
\[
g_n(x) = \phi_n^*(x : x_1,\ldots,x_n) - f(x) 
       = \hat{h}_{n,k_n}(x) - f(x) ,
\]
and note that for all $1 \leq i \leq k_n$,  
\be
\label{gnlessfivealpha}
V(g_n:-i,i) \, \leq \, V(\phi_n^*:-i,i) + V(f:-i,i) 
\, < \, 5 \alpha(i). 
\ee
Fix $\epsilon > 0$.  Select an integer $L \geq 1$ such that 
\be
\label{tail}
\int_{|x| \geq L} f(x) dx \leq \epsilon 
\ee
and define
\be
\label{delta}
\delta = \frac{\epsilon}{L} \, .
\ee
Finally, choose an integer $K \geq 1$ so large that 
\be
\label{interval}
2^{-K} < \frac{\epsilon \delta}{\alpha(L) (50 \alpha(L) + 5 \delta)}.
\ee

\vskip.07in

As $\bx \in \Omega(f)$ and the partitions $\pi_k$ are nested, 
there exists an integer $N = N(\bx, \epsilon, f, \alpha)$
such that for $n \geq N$ one has $k_n \geq \max\{K,L\}$,
\be
\label{intgn}
|\int_A g_n(x) dx| \, = \, |\hat{\mu}_n(A) - \mu_f(A)| 
                   \, < \, \frac{\delta}{2} \cdot 2^{-K}
\ee
for $A \in \pi_K$ with $A \subseteq [-L,L)$, and
\be
\label{approx}
|\hat{\mu}_n\{ |x| \geq L \} - \mu\{ |x| \geq L \} | \leq \epsilon \, .
\ee
For each $n$ let
\[
H_n = \{ x \in \real : |g_n(x)| > \delta \}
\]
contain those points having large error, and let
\[
\label{defcalH}
\HC_n = \{ A \in \pi_K : A \cap H_n \neq \emptyset, A \subseteq [-L,L) \}.
\]

Fix $n \geq N$ and consider a set $A \in \HC_n$.
By definition, there exists a point
$x \in A$ such that $|g_n(x)| > \delta$.  Assume for the moment
that $g_n(x) > \delta$. 
It follows from (\ref{intgn}) that there is a
point $y \in A$ such that $g_n(y) < \delta/2$, and therefore
\be
\label{jump}
\sup_{x,y \in A} |g_n(x)-g_n(y)| > \delta / 2 \, .
\ee 
As $k_n \geq L$ the variation of 
$g_n$ on $A$ is less than $5 \alpha(L)$ by (\ref{gnlessfivealpha}),
so that for each $z \in A$,
\[
g_n(z) \leq g_n(y) + 5 \alpha(L) \leq \frac{\delta}{2} + 5 \alpha(L)
\, ,
\]
and
\[
g_n(z) \geq g_n(x) - 5 \alpha(L) \geq \frac{\delta}{2} - 5 \alpha(L)
\, .
\]
Therefore,
\be
\label{bound}
\sup_{z \in A} |g_n(z)| \leq \frac{\delta}{2} + 5 \alpha(L) \, .
\ee
A similar argument in the case $g_n(x) < - \delta$ shows that
both (\ref{jump}) and (\ref{bound}) hold for each $A \in \HC_n$.
It is immediate from (\ref{jump}) that
\[
\frac{\delta}{2} |{\cal H}_n| \leq V(g_n: -L,L) < 5 \alpha(L) ,
\]
and consequently
\be
\label{number}
 |{\cal H}_n| < \frac{10 \alpha(L)}{\delta}.
\ee

\vskip.06in

For each $n \geq N$ the integrated error between $\phi_n^*$ and $f$ may be
decomposed as follows:
\begin{eqnarray*}
\lefteqn{\int |\phi_n^*(x) - f(x)| dx}  \\
& \leq & \sum_{A \in {\cal H}_n} \int_{A} |g_n(x)| dx
  \  +   \sum_{A \notin {\cal H}_n, A \subseteq [-L,L)} \int_{A} |g_n(x)| dx 
  \  + \, \int_{|x| \geq L} |g_n(x)| dx \\
& \deq & \Theta_1 + \Theta_2 + \Theta_3
\end{eqnarray*}
Inequalities (\ref{bound}), (\ref{number}) and (\ref{interval})
imply that 
\[
\Theta_1 \, \leq \,
\sum_{A \in {\cal H}_n} \int_{A} (\frac{\delta}{2} + 5 \alpha(L)) dx
\, \leq \,
\left( 5 \alpha(L) + \frac{\delta}{2} \right) 
\frac{10 \alpha(L)}{\delta 2^K}
\, \leq \, \epsilon \, ,
\]			
and by virtue of (\ref{delta}),
\[
\Theta_2 \, \leq \, 
 \int_{[-L,L)} \delta dx
= \delta \cdot 2 L = 2 \epsilon \, .
\]
Finally, it follows from (\ref{tail}) and 
(\ref{approx}) that 
\begin{eqnarray*}
\Theta_3 & \leq &
\hat{\mu}_n \{ |x| \geq L \} + \mu \{|x| \geq L \}
\leq 3\epsilon.
\end{eqnarray*}
Combining these three bounds shows that 
\[
\limsup_{n \to \infty} \int |\phi_n^*(x) - f(x))| dx \leq 6 \epsilon \, ,
\]
and as $\epsilon$ was arbitrary, the desired $L_1$ convergence of 
$\phi_n^*$ to $f$ follows.  $\Box$.

\section{Proof of Theorem \protect{\ref{thm2}}}
\label{PT2}

The following result can be established by a straightforward extension
of the Glivenko Cantelli Theorem, or by a bracketing argument
(c.f.\ Pollard \cite{Pol84}).

\begin{lem}
\label{GC}
Let ${\cal A}$ be the collection of all (finite and infinite)
intervals in $\real$.
If $\bx \in \Omega(f)$ then 
\[
\sup_{A \in {\cal A}} | \hat{\mu}_n(A) - \mu_f(A) | \, \to \, 0 \, .
\]
\end{lem}

\vskip.15in

\noindent
{\bf Proof of Theorem \ref{thm2}:}
Consider the family $\F_0 = \{ h_1, h_2, \ldots \} \subseteq \F$
of Rademacher densities where
\[
h_k(x) = 
\left\{ \begin{array}{ll} 
        2 & \mbox{ if $2j 2^{-k} \leq x < (2j+1) 2^{-k}$ for some
                   $0 \leq j < 2^{k-1}$ } \\ 
        0  & \mbox{otherwise .} \end{array} \right. 
\]
Note that each $h_j$ is supported on $[0,1]$ and that 
$\int |h_j(x) - h_k(x)| dx = 1$ whenever $j \neq k$.
Let $\mu_k$ be the probability measure having density $h_k$, and
for each finite sequence $u_1, \ldots, u_m \in [0,1]$ let
\[
\Delta_k(u_1,\ldots,u_m) = 
\sup_{A \in {\cal A}}
\left| \frac{1}{m} \sum_{j=1}^{m} I_A(u_j) - \mu_k(A) \right| \, ,
\]
measure the distance between $\mu_k$ and 
the empirical measure of $u_1, \ldots, u_m$.

We show that if $\Phi$ is 
consistent for $\F_0$ then there is a stationary sequence $\bx^*$
whose limiting density is identically one on $[0,1]$, but is such that
$\phi(\cdot : x_1^*, \ldots, x_n^*)$ fails to have a limit in $L_1$.
For each $k \geq 1$ select a sequence 
$\bx^{(k)} = (x_1^{(k)}, x_2^{(k)}, \ldots) \in \Omega(h_k)$
(e.g.\ a typical sample sequence from an i.i.d.\ process with
density $h_k$), and define 
\[
m_k = \min\left\{ M : 
\sup_{m \geq M} 
\Delta_k(x_1^{(k)}, \ldots, x_m^{(k)})
\leq \frac{1}{k+1}  \right\} \, .
\]
Lemma \ref{GC} insures that $m_k$ exists and is finite.

Fix any procedure $\Phi = \{\phi_1, \phi_2, \ldots\}$ 
that is consistent for ${\cal F}_0$ and
consider the infinite sequence $\bx^{(1)}$.  
As $h_1 \in {\cal F}_0$,
\[
\int | \phi_n(x : x_1^{(1)}, \ldots, x_n^{(1)}) - h_1(x) | dx \to 0
\]
as $n \to \infty$.
Therefore there is an integer
$n_1 \geq m_2$ and a corresponding initial segment
$\by^{(1)} = x_1^{(1)}, \ldots, x_{n_1}^{(1)}$ of $\bx^{(1)}$ 
such that 
\[
\int | \phi_{n_1}(x : \by^{(1)}) - h_1(x) | dx \leq \frac{1}{4}
\ \ \mbox{ and } \ \ \Delta_1(\by^{(1)}) \leq \frac{1}{2} \, .
\]

Now suppose that one has constructed a sequence $\by^{(k)}$ of 
finite length $n_{k}$ from initial segments of $\bx^{(1)}, \ldots, \bx^{(k)}$
such that 
\be
\label{close}
\int | \phi_{n_k}(x : \by^{(k)}) - h_k(x) | dx \leq 1/4 \, ,
\ee
\be
\label{supint}
\Delta_k(\by^{(k)}) \leq (k+1)^{-1} ,
\ee
and
\be
\label{nk}
n_k \geq k \cdot m_{k+1} \, .
\ee
As $\by^{(k)}$ is finite, the concatenation 
$\by^{(k)} \cdot \bx^{(k+1)}$ is contained in
$\Omega(h_{k+1})$.  It follows from the
consistency of $\Phi$ and Lemma \ref{GC} that when $n$ 
is large enough each initial segment 
$\by^{(k+1)} = \by^{(k)} \cdot (x_1^{(k+1)}, \ldots, x_{n-n_k}^{(k+1)})$
of $\by^{(k)} \cdot \bx^{k+1}$ satisfies 
(\ref{close}) and (\ref{supint}) with $k$ replaced by $k+1$.
Select $n_{k+1} > n_k$ so large that the same is true of (\ref{nk}).

As $\by^{(k+1)}$ is a proper extension of $\by^{(k)}$, repeating the
above process indefinitely yields an infinite sequence
$\bx^*$.  By construction, the functions
$\phi_n (\cdot) = \phi(\cdot : x_1^*, \ldots, x_n^*)$ do not
converge in $L_1$.  Indeed, it follows from (\ref{close}) 
and the triangle inequality that 
$\int |\phi_{n_k} - \phi_{n_l}| dx \geq 1/2$ whenever $k \neq l$.

It remains to show that the limiting density of $\bx^*$  
is uniform on $[0,1]$.  To this end, 
fix $k \geq 1$ and let $A \subseteq [0,1]$ 
be an interval of length $l(A)$.
It is easily verified that 
\be
\label{mukclosetolby1overk}
|\mu_k(A) - l(A)| \leq 2^{-k+1} \le \frac{1}{k} \, .
\ee
Let $\hat{\mu}_n(A)$ be the empirical distribution
of $A$ under $x_1^*, \ldots, x_n^*$, and for each 
$1 \leq r \leq n_{k+1} - n_k$ define
\[
\hat{\mu}_{r,k}'(A) = \frac{1}{r} \sum_{j=n_k+1}^{n_k+r} I_A(x_i^*)
\]
It follows from the equation
\[
\hat{\mu}_{n_k+r}(A) = 
\frac{n_k}{n_k+r} \cdot \hat{\mu}_{n_k}(A) + 
\frac{r}{n_k+r} \cdot \hat{\mu}_{r,k}'(A) 
\]
that the difference
\begin{eqnarray*}
\label{rat}
|\hat{\mu}_{n_k+r}(A) - l(A)| & \leq &
\frac{n_k}{n_k+r} \cdot |\hat{\mu}_{n_k}(A) - l(A)| \ + \ 
\frac{r}{n_k+r} \cdot |\hat{\mu}_{r,k}'(A) - l(A)| \\
& \deq & I + II .
\end{eqnarray*}
By virtue of (\ref{supint}) and (\ref{mukclosetolby1overk}),
\[
I \, \leq \,
|\hat{\mu}_{n_k}(A) - \mu_{k}(A)| + |l(A) - \mu_{k}(A)| 
\, \leq \, \frac{1}{k+1} + \frac{1}{k} \, .
\]
If $n_{k+1}-n_k\geq r \geq m_{k+1}$ then 
\[
\Delta_{k+1}(x_{n_k+1}^*, \ldots, x_{n_k+r}^*)
\, = \, \Delta_{k+1}(x_1^{(k+1)}, \ldots, x_r^{(k+1)})
\, \leq \, \frac{1}{k+2}
\]
and therefore
\[
II \, \leq \,
|\hat{\mu}_{r,k}'(A) - \mu_{k+1}(A)| \, + \,
|\mu_{k+1}(A) - l(A)| \, \leq \,
\frac{1}{k+2} + \frac{1}{k+1} .
\]
On the other hand, if $1\le r < m_{k+1}$ then (\ref{nk}) implies that
\[
II \, \leq \, \frac{2r}{n_k + r} 
\, \leq \, \frac{2r}{kr + r}
\, = \, \frac{2}{k+1} \, .
\]
These bounds insure that 
\[
\max\{ |\hat{\mu}_n(A) - l(A)| : n_k < n \leq n_{k+1} \}
\leq \frac{4}{k} ,
\]
and consequently
\[
\lim_{n \to \infty} |\hat{\mu}_n(A) - l(A)| = 0 \, .
\]
As $A \in {\cal A}$ was arbitrary, $\bx^*$ is stationary with limiting
density $f(x) = 1$ on $[0,1]$.  $\Box$

\section*{Acknowledgments}

The authors wish to thank L\'{a}szl\'{o} Gy\"{o}rfi 
for his helpful comments and
suggestions.



\small{

}

\end{document}